\documentclass{article}
\usepackage[utf8]{inputenc}
\usepackage[a4paper]{geometry}
\usepackage{amsmath}
\usepackage{amssymb}
\usepackage{amsthm}
\usepackage{graphicx}
\usepackage{hyperref}
\usepackage{color}
\usepackage{enumitem}
\theoremstyle{definition}
\newtheorem{definition}{Definition}[section]
\newtheorem{theorem}{Theorem}[section]

\newtheorem{proposition}{Proposition}[section]
\newtheorem{corollary}{Corollary}[section]
\title{The limit of human intelligence }
\author{ Santanu Acharjee$^{1, \dagger}$ and Upashana Gogoi$^{2,3}$\\
$^{1,2}$ Department of Mathematics\\
Gauhati University\\
Guwahati-781014, Assam, India\\
$^{3}$ Department of Mathematics\\
Morigaon College\\
Morigaon-782105, Assam, India\\
e-mail: $^1$sacharjee326@gmail.com,\\ $^2$upashanagogoi122@gmail.com\\
$^\dagger$Corresponding author: Santanu Acharjee}
\date{}
\begin{document}
\maketitle
{\bf Abstract:}  In 1998, Fields medalist Stephen Smale [S. Smale, Mathematical problems for the next century,  The mathematical Intelligencer, 20(2) (1998), 7-15] proposed his famous eighteen problems to the mathematicians of this century. The statement of his eighteenth problem is very simple but very important. He asked ``What are the limits of intelligence, both artificial and human?". In this paper, we prove  that  human intelligence is limitless.  Moreover, we provide justifications to state that  artificial intelligence has limitations. Thus, human intelligence will always remain superior to artificial intelligence. Moreover, we provide justifications to conclude the limitations of artificial intelligence. \\

{\bf 2020 AMS Classifications:} 91E10, 68T01, 68T27, 54E99.\\


\section{Introduction}
In \cite{1}, Smale proposed a list of eighteen important mathematical problems for the next century. There he listed down the eighteenth problem as ``What are the limits of intelligence, both artificial and human?". The question is simple but significant in light of recent developments in artificial intelligence (AI). The discussions surrounding the potential or  restrictions of artificial intelligence have become crucial study areas with its advancement. Certainly, AI possesses some abilities of the human mind, and it is an obvious curiosity of men to know about the extent of AI's potential. Regarding the limitations of AI, Wang \cite{2} discussed that there are three prevalent misunderstandings, namely, thinking of AI as having an axiomatic system, an AI system's approach to problem-solving  being  comparable to a Turing machine, or taking an AI system to be formal. Kelly \cite{3} mentioned that the future of AI is going to be cognitive, not ``artificial". Cognitive computing \cite{3} is the system that learns through its interactions and experiences with people and the environment rather than being explicitly designed. After the tabulating era and the programming era, the cognitive era \cite{3} began with the introduction of Watson, IBM's cognitive computing system, in 2011, which defeated Ken Jennings and Brad Rutter at the American television game show `Jeopardy!'. The ability to make sense of unstructured data is what makes cognitive computing the future of computing.\\ 

The foundation of psychology is made up of ideas from common sense understandings of mental and behavioural processes \cite{4}. It is more difficult to define psychological constructs than physical, biological, or chemical ones. Smedslund \cite{5} came to the conclusion that psychology cannot be an empirical science because of the irreversibility of psychological processes, individual variability of brains, which makes comparison and stability of results quite challenging, the infinite possibilities of people perceiving the world, etc. What makes psychology interesting is the shift towards behaviourism rather than structuralism. Lewin \cite{4} described psychological processes in topological spaces, concentrating on conceptualizations and quantifications of psychological forces. Lewin developed field theory to explain behaviour as the outcome of interactions between the individual and the environment. According to his theory, cognitive dynamics may be thought of as ``life space", which is a field that encompasses a person's values, wants, objectives, reasons, emotions, concerns, etc. He defined a functional relationship $B=f(P, E)$ which gives the resultant behaviour when a person $P$ and an environment $E$ interact, where both $P$ and $E$ belong to a psychological-topological space. Thagard \cite{6} stated, ``Most cognitive scientists agree that knowledge in the mind consists of mental representations." He further stated, ``Cognitive science proposes that people have mental procedures that operate on mental representations to produce thought and action." Eysenck and Keane \cite{7} considered the psychological space as a space of mental representations of some aspects of things from the physical world or things imagined. Sims et al. \cite{8} studied the structures and dynamics of psychological phenomena via topology. They described the cognitive space as a generalised algebraic mental structure $\mathcal{C}=(C, \sigma, I)$, where $C$ is the set of symbols or mental representations, the signature $\sigma$ describes some concatenation rules on $C$, and $I$ is an interpretation semantic function that gives meaning to representations and concatenation. In general, $C$ is called a cognitive or psychological space. \\

In 1936, Tarski \cite{9} drew up a formal formulation of the idea of logical consequence, which, according to him, matched quite well with the common understanding of consequence. This theory uses an operator $Cn$, called the consequence operator, which is a function defined on the power set of a given set of objects. Tarski referred to these objects as ``meaningful sentences" following a suggestion of his mentor Les\'{n}iewski \cite{10}. It is possible to produce certain other sentences from any set $A$ of sentences by using rules of inference. The collection of these sentences is called the consequence of the set $A$ denoted by $Cn(A)$. Tarski introduced a wide theory, as he did not specify the nature of these sentences; they can be in any type of scientific language \cite{10}. Logic is very closely related to cognitive processes \cite{11}. But mathematical logic encounters a number of difficulties when attempting to describe cognition and intelligence \cite{12}. Mathematical logic and cognitive logic are essentially different in many aspects. In mathematical logic, a term's interpretation determines its meaning; therefore, its meaning remains constant throughout the operation \cite{12}. On the contrary, a term's meaning in the human mind frequently varies depending on the experiences and contexts \cite{12}. Unlike mathematics, a compound term's meaning in the human mind cannot always be reduced to that of its constituent parts \cite{12}. For instance, although the word `keyword' consists of  two different words, `key' and `word', it never gives the picture of `key' and `word' independently in someone's mind, as depicted in figure \ref{fig2}. The meaning of the term `keyword' is entirely distinct from the meanings of the words that make up the word. A statement is either true or false in mathematical logic, but in the real world, the human mind cannot decide the truthiness of certain statements precisely \cite{12}. Moreover, people also change their minds about a statement's truth value after learning new facts about it, while in mathematical logic the truth value of the statement remains constant \cite{12}. In mathematical logic, inference processes adhere to some predefined algorithms, but human thinking does not always stick to a particular algorithm, and therefore inferences can be unpredictable and can be headed in unforeseen directions \cite{12}. Moreover, in traditional reasoning, a precise explanation of a result exists, but the human mind often comes to conclusions solely based on intuition and guesses \cite{12}. Traditional logic's inference rules are deductive in nature \cite{12}. In a deductive system, the information for the conclusion is already present in the premises. The inference rules make it clear. But in the case of human reasoning, it is not always a deductive system. There are some situations where the conclusion cannot be deduced from the premises alone. When there is not enough information or tools to use for deduction, intelligence is required \cite{12}. To deal with the limitations of traditional axiomatic systems, Wang \cite{13} designed a non-axiomatic reasoning system so that it could be adaptive and function even when there is a lack of information and resources.\\ 

Woleński \cite{11} divided a theory $T$ into its starting assumptions and their implications, and  he took the theory $T$ as an organised collection of sentences. This division and the process of inferring the consequences from the premises are involved with the operation of logical consequence. If $X \subseteq T$ is a collection of initial assumptions, then $T=Cn(X)$ means $T$ is a set of logical consequences of X. From a metamathematical point of view, a theory is a set of sentences that is closed by the consequence operation, i.e., $Cn(T)\subseteq T$. Moreover, this consequence operation $Cn$ satisfies Tarski's general axioms \cite{14} from which $T\subseteq Cn(T)$ follows. Thus, $T=Cn(T)$, i.e., a theory is a set of sentences that is equal to the set of its logical consequences. Woleński \cite{11} further defined logic as a theory. Logic generates conclusions from some initial assumptions. The statement ``$B$ is provable from the assumptions $A$" is formally written as $A \vdash B$. He assumed that the consequence operation satisfied the deduction theorem: if $B\in Cn(X \cup \{A\})$, then $(A\Rightarrow B) \in Cn(X)$. From this result, $A \vdash B$ gives $\vdash (A \Rightarrow B)$. Which can be written as $\phi \vdash (A \Rightarrow B)$. That means the formula $A\Rightarrow B$ in logic can be derived from an empty set of assumptions. This motivates the metalogical definition of logic as a theory, formally written as $\mathbf{LOG}=Cn(\phi)$. So, logic is the consequence of an empty set of assumptions. In order to prove a logical rule, we do not need something logical in the first place as a presupposition. Here lies the importance of the theory of logical consequence in the study of cognitive science. Recently, Muravitsky \cite{15} studied monotonic reasoning by using  consequence relations. In this paper, we introduce a mental structure $\mathcal{C}=(C, \sigma, I, Cn)$ called the cognitive-consequence space. Moreover, we construct a cognitive-consequence topological space on $\mathcal{C}$ and study some fundamental properties, including cognitive closure and convergence of sequences of thoughts in $\mathcal{C}$. We introduce a cognitive filter and a cognitive ideal and construct a new topological space from the perspective of cognition. We also introduce the concept of G\"{o}del's incompleteness black hole in the solution space of a problem inside $\mathcal{C}$.\\

\section{Preliminaries}
In this section, we discuss some preliminary definitions from \cite{8}. These definitions are defined on the mental substructure $\mathcal{C}=(C, \sigma, I)$ and the dynamic 4-dimensional physical time-space world $R^3(t)$.\\

\begin{definition} \cite{8}
The neighborhood $N(x)$ of an element $x$ in $R^3(t)$ is defined to be ``a set of parts as a practical working whole" that contains $x$. Here, the set $N$ of parts in $R^3(t)$ is considered to be a neighborhood of any of its “parts” if and only if that set $N$ as a whole satisfies a practical need of a person.
\end{definition}

\begin{definition}\cite{8}
In some environment $E$ of $R^3(t)$, let B be a collection of objects $x \in E$, where $x$ is either complete but not a practical whole, or an irreducible set that is a practical whole. B is a base for a neighborhood system in $E$ if and only if every practical whole $K$ in $E$ contains some object(s) $x$ in B.
\end{definition}

\begin{definition} \cite{8}
A complete object in $R^3(t)$ is an object that is considered to be whole, a unit, with respect to some psychology.
\end{definition}

\begin{definition} \cite{8}
An irreducible set is a connected set that is not the union of any other complete sets.
\end{definition}

\begin{definition} \cite{8}
Let B be a base for a neighborhood system. A set $K$ is open if and only if $K$ is a finite union of objects from B and $K$ is a practical whole.
\end{definition}

\begin{definition} \cite{8}
Let E be some environment in $R^3(t)$. Let B=$\{B_i|\  B_i$ is complete or a connected practical whole in $E\}$. Let $K_j=\bigcup B_i$, for $B_i\in L$, where $L$ is some finite subset of B, and $K_j$ is a practical whole in $E$. Then, the collection of sets $T=\{\phi, E\}\bigcup \{K_j\}$ is a topology on $E$ and B is a base for $T$.
\end{definition}

\begin{definition}
\cite{8} The closure of any base object $\{B_n\}$ is defined to be the practical-whole union, of which $B_n$ is a part of; that is, $\overline{\{B_n\}}=K_j=\bigcup B_i$, where $K_j$ is a practical whole. In this sense, each $B_i$ is a “limit object” of $K_j$. Here, limit objects are defined to be those objects necessary in a union to make that union a practical whole (pw).
\end{definition}

\begin{proposition} \label{p2.1}
\cite{8} From definitions 2.5 and 2.7, the practical wholes, $K_j=\bigcup B_i$ in the topology defined in 2.6 are both open and closed (clopen).
\end{proposition}

\section{Mental space as a cognitive substructure}

Sims et al. \cite{8} expressed the mind-space by a mathematical structure, which was termed a mental substructure ($\mathcal{C}$). We cannot perceive a boundary between the physical world and the mind-space, as there are no neighbourhoods of mental points and physical things that intersect both of them. Further, it is not possible to define any kind of physical distance between the thoughts \cite{8}. Also, the thoughts do not possess any geometric shape. Thoughts are the basic building blocks of the mind-space, which can be something from the perceptual physical world in $R^3(t)$ or the mental representations of the things imagined in the mind itself \cite{7}. Thoughts can be combined to produce other thoughts, or new thoughts can be inferred from the preceding thoughts. While doing that, the concept of consequence theory plays an important role. Thoughts can be influenced by the consequences of the preceding thoughts. We further generalise this mental substructure by associating the consequence operator. We call it a cognitive-consequence space. Mathematically, we can denote the cognitive-consequence space as $\mathcal{C}=(C, \sigma, I, Cn)$, where $C$ is the set of all symbols, or mental representations of things from $R^3(t)$ or the mind itself, $\sigma$ describes some mental-grammar concatenation rules on $C$, $I$ interprets and gives meanings to the representations and concatenations, and $Cn$ is the consequence operator. Throughout this paper, we refer to any mental activity, event, formula, thought, etc. as a mental representation.\\

We procure  the definition of consequence operator from \cite{11}. The consequence operator is  a function $Cn$:$\ P(C)\rightarrow P(C)$ that satisfies the following properties:
\begin{enumerate}[label=(\roman*)]
    \item denumerability of the language as a set of sentences,
    \item  $A \subseteq Cn(A)$  (the inclusion axiom),
    \item  $A \subseteq B \Rightarrow Cn(A) \subseteq Cn(B)$  (monotonicity),
    \item  $Cn(Cn(A))= Cn(A)$ (idempotence of $Cn$),
    \item if $X \in Cn(A)$; then there exists a finite set $B \subseteq A$ such that $X \in Cn(B)$,
    \item if $Y \in Cn(A \cup \{X\})$; then $(X \Rightarrow Y) \in Cn(A)$ (deduction theorem),
\end{enumerate}
for any two subsets $A$ and $B$ of $C$. Here, $P(C)$ indicates the power set of $C$.\\ 

The above axioms are due to Tarski \cite{16}. We can show the denumerability of the thoughts present in the cognitive-consequence space at a given instant of time. For instance, consider the situation $f_1: ``x_1$ {\it is playing} $x_2$" in $\mathcal{C}$. Here, $x_1$ can be a mental representation of a person or animal in $R^3(t)$, $x_2$ can be a mental representation of a game in the same space, and ``is playing" is also a mental representation of a task in $R^3(t)$. Taking $x_1$, $x_2$, and ``is playing" to be the basic elements in the situation, i.e., $x_1$, $x_2$, ``is playing" $\in C$, we get that the situation $f_1 \in C$. In other words, a situation can be any mental representation formed in the mind having a syntax $\sigma$ and interpretation $I$ on it. We can also combine two or more mental representations to form compound mental representations, and new mental representations can also be inferred from the old ones. Thus, the subsets of $C$ are collections of some mental representations from the cognitive-consequence space. Now, we define the following definitions and some related examples:\\ 

\begin{definition}
 In a cognitive-consequence space $\mathcal{C}=(C, \sigma, I, Cn)$, a  subset $A$ of  $C$ is said to be deductive if $Cn(A)=A$.
\end{definition}

\begin{definition}
 Consider $\mathcal{C}=(C,\sigma, I, Cn)$ be a cognitive-consequence space. Then, the collection $\tau = \{ A \subseteq C: Cn(C-A)=(C-A)\}$  is  called a cognitive-consequence topology (in short, $CCT$) on $C$. Moreover, we call  $(C,\tau)$ as a cognitive-consequence topological space (in short, $CCTS$).\\

The cognitive space $C$ is not a deductive system in nature \cite{17}. Thus, $\phi \notin \tau$ because $Cn(C-\phi)=Cn(C)\neq C$.
Also, $C\notin \tau$, as $Cn(C-C)=Cn(\phi)\neq\phi$, i.e., the consequence of the empty set is not empty, unlike topological closure \cite{8}.\\

\end{definition}

\begin{definition}
 A subset $A$ of $C$ in a  $CCTS$ $(C,\tau)$ is said to be  consequence-wise open (CWO) if $A\in \tau$. On the other hand, $A^c$ is said to be consequence-wise closed (CWC).
\end{definition}

\begin{theorem} \label{t1}
In $(C, \tau)$, the arbitrary union of CWO sets is CWO.
\end{theorem}
\begin{proof}
Consider $\{A_i|\; i \in \Delta, \Delta \;$ is an index set$\} \subseteq \tau$. Then, $Cn(C-A_i)=C-A_i \forall i\in \Delta$. \\
Now,

\begin{equation*} 
\begin{split}
  Cn(C-(\bigcup_{i \in \Delta} A_i)) & =Cn(\bigcap_{i \in \Delta} (C-A_i))\\
  & = \bigcap_{i \in \Delta} Cn(C-A_i)\\
  & = \bigcap_{i \in \Delta} (C-A_i) \\
  & = C-(\bigcup_{i  \in \Delta} A_i)    
\end{split}
\end{equation*}
Thus,  $\bigcup \limits_{i \in \Delta} A_i\in \tau$. Hence, theorem is proved.
\end{proof}

\noindent In general, in a cognitive-consequence space, the property $Cn(A\bigcup B)=Cn(A)\bigcup Cn(B)$ does not hold \cite{16}. The equality holds if $A$, $B$ and $(A \bigcup B)$ are deductive systems, as if it happens, then we get $Cn(A\bigcup B) = A \bigcup B=Cn(A) \bigcup Cn(B)$. In a cognitive sense, if we consider two subsets of mental representations in $C$, where the consequences of these two sets already exist inside the sets, then the consequence we obtain by combining these two sets will be equal to the union of the individual consequences. For example, consider a set $S=$\{l, r\} of pair of socks where l denotes the left sock and r denotes the right sock, and let $A$=\{wear the sock l\} and $B=$ \{wear the sock r\} be two subsets of the cognitive space $C$, then $A$ and $B$ will form deductive systems individually. Combining $A$ and $B$ together, we get $A\bigcup B=$\{wear the sock l, wear the sock r\}, which  will again form a deductive system, and the consequence of $A \bigcup B$ will be equal to the union of the consequences of $A$ and $B$. In \cite{8}, different examples were taken in order to explain the concept of practical working whole. In a $CCTS$, we are considering those subsets of $C$ whose complements are deductive systems. If we consider an arbitrary collection of subsets of a $CCTS$, then their intersection may not be present in the $CCTS$ unless the complement of the intersection forms a deductive system,  which is nothing but the arbitrary union of some deductive systems.    \\ 

\begin{theorem} \label{t2}
In $(C, \tau)$, the arbitrary intersection of CWO sets is CWO if the arbitrary union of the complements of the CWO sets forms a deductive system. 
\end{theorem}
\begin{proof}
Consider $\{A_i|\; i \in \Delta, \Delta$ is an index set$\} \subseteq \tau$ be an arbitrary collection of CWO sets, then $Cn(C-A_i)=(C-A_i) \forall i\in \Delta$. Now, $Cn(C-(\bigcap \limits_{i \in \Delta} A_i))=Cn(\bigcup \limits_{i \in \Delta}(C-A_i))$. Since, $\bigcup \limits_{i \in \Delta}(C-A_i)$  is a deductive system, hence we have $Cn(\bigcup \limits_{i \in \Delta}(C-A_i))= \bigcup \limits_{i \in \Delta}(C-A_i)$.\\
So, $Cn(C-(\bigcap \limits_{i \in \Delta} A_i))=Cn(\bigcup \limits_{i \in \Delta}(C-A_i))= \bigcup \limits_{i \in \Delta}(C-A_i)=C-\bigcap\limits_{i \in \Delta} A_i$. Thus, $\bigcap\limits_{i \in \Delta} A_i$ is a CWO set.
\end{proof}

\noindent From the above results, it is well understood that $\tau$ is neither a topology nor an Alexandrov topology on $C$. Thus, $(C,\tau)$ is expected to provide different properties than those of the properties of general topological space or Alexandrov topological space. The human mind cannot always make deductive inferences from a given set of mental representations in the cognitive space \cite{17}. Also, the thoughts and mental representations depend on time and the environment. New thoughts enter and leave the human mind with the passing of time, which makes the mind an open system \cite{2}. On the other hand, there may exist some subsets of $C$ that form deductive systems. Now, what are the thoughts that constitute the $CCTS$ is an important question to ask. \\

\begin{theorem} \label{t3}
In $(C,\tau)$, there always exists at least one mental representation $f\in C$ that does not belong to any CWO set. 
\end{theorem}
\begin{proof}
If every mental representation $f$ belongs to some CWO sets in $\tau$, by theorem \ref{t1},  we  get $C\in \tau$, which is a contradiction.
\end{proof}

In \cite{4}, Lewin urged  the importance of overlapping  two regions in the psychological space. He stated, ``{\it  A psychologically important application of the concept of
overlapping is the overlapping of two situations. A child may be eating and, at the same time, listening to the song of a bird.
The listening can be the major and the eating the minor activity, or the reverse. Between the two extremes many transitions are possible. Such cases in which one is involved to different
degrees in two different activities are of common occurrence. But they offer considerable difficulty for description as well as for treatment of their dynamic facts}". Thus, it is clearly understood that any two CWO sets are not disjoint in  anyway. Thus, chances of overlapping of two CWO sets are always present.  Thus, we obtain the following theorem:\\

\begin{theorem} \label{t4}
In $(C, \tau)$, there always exists at least one mental representation $f \in C$ that belongs to a CWO set.
\end{theorem}
\begin{proof}
Let $f \in C$ and no $f$ belongs to any CWO set $A_i$ $ \forall i \in \Delta$. \\
Then, \begin{eqnarray*}
    f \notin A_i & \Rightarrow & f \in C-A_i\\
    & \Rightarrow & \bigcup_{f \in C} \{f\} \subseteq C - \bigcap_{i \in \Delta} A_i\\
    & \Rightarrow & C \subseteq C - \bigcap_{i \in \Delta} A_i \\
    & \Rightarrow & \bigcap_{i \in \Delta} A_i = \phi
    \end{eqnarray*}
which is a contradiction. 
\end{proof}

\begin{corollary}
The intersection of all the deductive systems whose complements form a $CCTS$ is non-empty.
\end{corollary}
\begin{proof}
The mental representation $f$ which does not belong to any of the sets $A$ in $\tau$ will be at the intersection of all the deductive systems $(C-A)$.
\end{proof}

\begin{theorem} \label{t5}
If a collection $A$ of mental representations and its complement in the cognitive-consequence space both form deductive systems, then the set $A$ cannot belong to the $CCT$. 
\end{theorem}
\begin{proof}
Let $A$ and $(C-A)$ both form deductive systems in the cognitive-consequence space $\mathcal{C}$, i.e., $Cn(A)=A$ and $Cn(C-A)=(C-A)$. From this, we get that both $A$ and $(C-A)$ belong to the $CCT$. Now, for every mental representation $f \in C$, we get either $f \in A$ or $f \in (C-A)$. Due to theorem \ref{t3}, we arrive at a contradiction.
\end{proof}

\begin{definition}
Let $(C,\tau)$ be a $CCTS$. Then, the cognitive closure of a set $A \subseteq C$ is the intersection of all the deductive systems that contain $A$. It is denoted by $Cl^\Box(A)$. 
\end{definition}

\begin{theorem}
Let $(C, \tau)$ be a $CCTS$ and $A \subseteq C$. Then, the cognitive closure $Cl^\Box(A)$ is the smallest deductive system that contains $A$.
\end{theorem}
\begin{proof}
Let $(C, \tau)$ be a $CCTS$ and $A \subseteq C$. Then, $Cl^\Box(A)=\bigcap\{A_i \subseteq C: Cn(A_i)=A_i, A \subseteq A_i, i \in \Delta \}$. Now, $Cn(\bigcap \limits_{i \in \Delta} A_i) \subseteq \bigcap \limits_{i \in \Delta} Cn(A_i) = \bigcap \limits_{i \in \Delta} A_i $. Thus, $\bigcap \limits_{i \in \Delta} A_i $ is a deductive system that contains $A$.\\

Now, we show that $Cl^\Box(A)=\bigcap \limits_{i \in \Delta} A_i$ is the smallest deductive system that contains $A$. For that, consider $B$ be the smallest deductive system such that $A \subseteq B \subseteq Cl^\Box(A)$. Since $Cn(B)=B$ and $A \subseteq B$, so $B=A_i$ for some $i \in \Delta$. Thus,  $\bigcap \limits_{i \in \Delta} A_i \subseteq B$ or $Cl^\Box(A) \subseteq B$. Thus, we get $B=Cl^\Box(A)$, i.e., $Cl^\Box(A)$ is the smallest deductive system that contains $A$.
\end{proof}

\begin{corollary}
The cognitive closure of the empty set is not empty.
\end{corollary} 
\begin{proof}
 The empty set is not a deductive system since $Cn(\phi) \neq \phi$ \cite{11}. 
\end{proof}
 \noindent From the above corollary, we can say that the cognitive closure is not as same as the topological closure. Also, monotonic consequence operator is not same as Kuratowski closure axioms.
 
\begin{theorem} \label{t6}
The consequence of a set of mental representations $A$ in a $CCTS$ $(C, \tau)$ is the smallest deductive system that contains $A$. 
\end{theorem}
\begin{proof}
We have, 
\begin{equation}
  A \subseteq Cn(A) \Rightarrow Cl^\Box(A) \subseteq Cl^\Box(Cn(A)) \Rightarrow Cl^\Box(A) \subseteq Cn(A). 
\end{equation}
since $Cn(A)$ is a deductive system and it is the smallest deductive system that contains itself, so,  $Cl^\Box(Cn(A))=Cn(A)$. \\
Again, 
\begin{equation}
  A \subseteq Cl^\Box(A) \Rightarrow Cn(A) \subseteq Cn(Cl^\Box(A)) \Rightarrow Cn(A) \subseteq Cl^\Box(A)
\end{equation}
since $Cl^\Box(A)$ is a deductive system, thus $Cn(Cl^\Box(A))=Cl^\Box(A)$.\\
From (1) and (2), $Cn(A)=Cl^\Box(A)$, which is the smallest deductive system that contains $A$.
\end{proof}

\begin{corollary} \label{c3}
The cognitive closure of a set $A \subseteq C$ of a $CCTS$ $(C, \tau)$ is equal to the consequence of the set, i.e. $Cn(A)=Cl^\Box(A)$. 
\end{corollary}
\begin{corollary}
$Cl^\Box(\phi)\neq \phi$ and $Cl^\Box(C)\neq C$.
\end{corollary}

\begin{corollary}
If $A \subseteq C$ of a $CCTS$ $(C, \tau)$ and $A\in \tau$, then $Cl^\Box(A) \neq A$ and $Cl^\Box(C-A)= (C-A)$.
\end{corollary}
\begin{proof}
Let $A \in \tau$. Then, by theorem \ref{t5} $Cn(A) \neq A$. Then, $Cn(C-A)=(C-A)$. From theorem \ref{t6}, $Cl^\Box(A)=Cn(A) \neq A$ and $Cl^\Box(C-A)=Cn(C-A)=(C-A)$.
\end{proof}

\noindent {\bf Properties:} For two sets $A, B \subset C$ of a cognitive-consequence topological space $(C, \tau)$, the following results hold: 
\begin{enumerate}[label=(\roman*)]
    \item if $A \subseteq B$, then $Cl^\Box(A) \subseteq Cl^\Box(B)$,
    \item $Cl^\Box(A) \bigcup Cl^\Box(B) \subseteq Cl^\Box(A \bigcup B)$,
    \item  $Cl^\Box(A \bigcup B)=Cl^\Box(A) \bigcup Cl^\Box(B)$ if $A, B$ and $(A \bigcup B)$ are deductive systems,
    \item $Cl^\Box(A \bigcap B) \subseteq Cl^\Box(A) \bigcap Cl^\Box(B)$,
    \item $Cl^\Box(A \bigcap B)=Cl^\Box(A) \bigcap Cl^\Box(B)$ if $A, B$ and $(A \bigcap B)$ are deductive systems.
\end{enumerate}

\begin{proof}
From corollary \ref{c3}, $Cl^\Box(A)=Cn(A)$. Hence, we get the following results:
\begin{enumerate}[label=(\roman*)]
    \item from property (iii) of the consequence operator if $A \subseteq B$, then $Cn(A) \subseteq Cn(B) \Rightarrow Cl^\Box(A) \subseteq Cl^\Box(B)$. 
    \item since $Cn(A) \bigcup Cn(B) \subseteq Cn(A \bigcup B)$. Thus, we get $Cl^\Box(A) \bigcup Cl^\Box(B) \subseteq Cl^\Box(A \bigcup B)$.
    \item Since $Cn(A) \bigcup Cn(B) = Cn(A \bigcup B)$ if $A, B$ and $(A \bigcup B)$ are deductive systems. So, we get $Cl^\Box(A) \bigcup Cl^\Box(B) = Cl^\Box(A \bigcup B)$ if $A, B$ and $(A \bigcup B)$ are deductive systems.
    \item Since $Cn(A \bigcap B) \subseteq Cn(A) \bigcap Cn(B)$. So, we get $Cl^\Box(A \bigcap B) \subseteq Cl^\Box(A) \bigcap Cl^\Box(B)$.
    \item Since $Cn(A \bigcap B)=Cn(A) \bigcap Cn(B)$ if $A, B$ and $(A \bigcap B)$ are deductive systems. So, we get $Cl^\Box(A \bigcap B)=Cl^\Box(A) \bigcap Cl^\Box(B)$ if $A, B$ and $(A \bigcap B)$ are deductive systems.
\end{enumerate}
\end{proof}

\subsection{Limits of thoughts}
In a cognitive-consequence space, a thought can lead to a new thought, which in turn leads to another thought, thus making a sequence of thoughts \cite{8}. A sequence of thoughts is different from the sequences defined in mathematics. A sequence of thoughts may converge to a single thought, just like a convergent sequence in mathematics. It may also diverge to different sequences of thoughts, resulting in different or the same conclusions at the same time by following different paths. It is also possible that there may be more than one sequence of thoughts that cognitively coincide to give a cognitive limit. Since the human thinking process does not always follow any specific rule \cite{12}, a sequence of thoughts can be headed in different directions depending on factors such as the state of the mind, the state of the surrounding environment, etc. The possibility that there may exist more than one cognitive limit of thoughts in the mind can be justified by the statement of Thagard \cite{6} about the functioning of the brain as {\it ``There is no single computational model of mind, since different kinds of computers and programming approaches suggest different ways in which the mind might work. The computers that most of us work with today are serial processors, performing one instruction at a time, but the brain and some recently developed computers are parallel processors, capable of doing many operations at once."} \\

In order to solve a problem or to make a decision, the human mind goes through a sequence of thoughts which ultimately leads to the solution. This process of arriving at a definite solution is not always a simple process. Sometimes a sequence of thoughts does not give a satisfied result or it becomes too difficult to proceed further. In that case, it is necessary to follow a different sequence of thoughts which may or may not give the desired result either. This process continues until the problem is solved or the decision is made. It may also possible that the desired result can be attained by following different sequence of thoughts. In that case, the question of feasibility comes into play. \\

\begin{definition}
Let $\mathcal{C}=(C,\sigma, I, Cn)$ be a cognitive-consequence space. The cognitive similarity distance on $C$ can be defined as a function $Cog$:$\ C \times C \rightarrow [0,1]$, which satisfies the following properties :
\begin{enumerate}[label=(\roman*)]
    \item $Cog(x,y) \geq 0$,
    \item $Cog(x,y)=0 \Leftrightarrow x \approx y$,
    \item $Cog(x,y)=Cog(y,x)$, 
    \item if $x \approx z$, then $Cog(x,y)=Cog(z,y)$,
    \item $Cog(x,z) \leq Cog(x,y)+Cog(y,z)$, where $x,y,z \in C$.
\end{enumerate}
Here, $x \approx y$ indicates that the two thoughts $x$ and $y$ are semantically similar \cite{18} or cognitively similar \cite{19} or identical \cite{19}. In short, we say that $x$ and $y$ cognitively coincide. The range is considered $[0,1]$ since most of the similarity measures are calculated between 0 and 1. For example, the Jaccard similarity measure \cite{20}, cosine similarity \cite{21}, etc. are calculated in the interval $[0,1]$. 
\end{definition}

\begin{definition}
A cognition ball $B(x,\epsilon)$ with a thought $x \in C$ as the centre and cognitive threshold $\epsilon \in (0,1)$ is the collection of all thoughts in $C$ whose cognitive similarity distance from $x$ is less than $\epsilon$, i.e., $B(x,\epsilon)=\{y \in C: Cog(x,y) < \epsilon\}$.
\end{definition}

\noindent The reason for taking the cognitive threshold $\epsilon$ to be in the open interval $(0,1)$ is that the threshold value depends on how we approach solving a problem rather than the solution itself. Moreover, whenever a new problem is given to someone, she tries to solve the problem with the help of some similar past experiences \cite{22}. This nature of problem solving from past experiences is not only applicable to human beings but also to other animals. Experiments show that chimpanzees \cite{23}, crows \cite{24}, etc. also use their past experiences while solving an unknown problem or dealing with any unfamiliar circumstances. Different types of mathematical measures like weight function \cite{25}, similarity measure \cite{26}, heuristic similarity measure \cite{27}, etc., and non-mathematical measures like cognitive similarity \cite{19}, psychological similarity \cite{28}, religious similarity \cite{29}, semantic similarity \cite{18}, shape similarity \cite{30}, etc. can be used to measure various cognitive distances. In the case of animals, birds, and insects, we have experimental evidence that they use various kinds of similarity measures to distinguish objects, shapes, etc.  For example, pigeons discriminate shapes by using topological similarity \cite{31}. Based on some parameters, such as the shape area and contour density, honey bees can distinguish different shapes \cite{32}. Also, experiments show that honey bees can discriminate shapes that are topologically different more rapidly than those that are topologically similar \cite{33}. The emotional similarity in the human mind depends on the co-occurrence of emotions happening in the day-to-day life \cite{34}. Thus, $Cog(x,y)$ indicates various types of cognitive similarity distances between two thoughts $x$ and $y$. Here, the thoughts may be mathematical or non-mathematical. Thus, $Cog(x,y)$ may be defined either mathematically or non-mathematically. For example, when we want to solve a mathematical problem, then our thoughts are purely mathematical. In this case, we may use  mathematical measures for cognitive similarity distance. But not all the thoughts in the human mind are mathematical, as there are emotional thoughts, cultural thoughts, religious thoughts, etc. \cite{17}. In these cases, we use non-mathematical measures in order to compute the cognitive similarity distances. Here, we consider the weight  as a particular case for computing the cognitive distance between two thoughts, which is a mathematical measure.\\

\begin{definition}
 The cognitive distance between two thoughts $x$ and $y$ in a cognitive-consequence space $\mathcal{C}=(C,\sigma, I, Cn)$ can be defined as $Cog(x,y)=| w(x)-w(y)|$, where $w$:$\ C \rightarrow [0,1]$ is a weight function, and $x,y \in C$.   Moreover, if $w(x)=w(y)$, then $x \approx y$.
\end{definition}
In a sequence of thoughts $\{x_i\}_{i=1}^{\infty}$, in a cognitive-consequence space $\mathcal{C}=(C,\sigma, I, Cn)$, we can associate a weight to every thought in the sequence that ultimately leads to the final thought, say $x$, where we consider $w(x)$ to be $1$. Initially, the sequence may be either finite or infinite, depending on the problem to be solved. The thoughts that are near the final thought $x$ will have greater weight than those further from $x$ in the sequence. For example, consider that we need to solve the mathematical equation $x^2-5x+6=0$. The steps $x_n$ to solve the problem have different weights in the sequence as shown below.
\begin{align*}
x_1& :  x^2-5x+6=0; \quad w(x_1)= 0.1,\\
x_2& :  x^2-2x-3x+6=0; \quad w(x_2)= 0.3,\\
x_3 &:  x(x-2)-3(x-2)=0; \quad w(x_3)=0.5,\\
x_4 &:  (x-2)(x-3)=0; \quad w(x_4)= 0.7,\\
x_5 &: {\text{either}} \; (x-2)=0 \; {\text{or}} \; (x-3)=0; \quad w(x_5)= 0.85,\\
x_6 &: {\text{either}} \; x=2 \; {\text{or}} \; x=3; \quad w(x_6)=0.9,\\
x_7 &: x=2,3; \quad w(x_7)=1.
\end{align*}
Here, $x_7$ is the final thought or solution to the problem, having a weight value of $1$. Note that the step $x_1$ can also be written as $x_1': 6-5x+ x^2=0$ which gives cognitively the same meaning as $x_1$ having the same weight value. Thus, $x_1 \approx x_1'$ since $x_1$ and $x_1'$ cognitively coincide with each other. Similarly, $x_2$ cognitively coincides with $x_2': 6-3x-2x+x^2=0$ or with any other representation giving the same cognitive meaning, which is written as $x_2 \approx x_2'$ and so on. It may be found that the above problem requires a finite number of steps to arrive at a solution. But, if we think about it more deeply, the process involves more than these few steps. To execute the first step, we first need to know what `$x$' means. For that, we need to know the language `English', and we must be familiar with the alphabets in the English language. Again, the letter $x$ acts as an unknown in the given equation. So, the concept of the unknown is also something to know. After that, we need to know many concepts, viz., digits',`square',`multiplication',`subtraction', etc. Although the solution of the mathematical equation shown in the above example appears to be solved in only six steps, the process actually involves infinitely many thoughts sequentially occurring within the mind within a short span of time. In general, when one wants to obtain a solution or conclusion on any matter, she tries to procure the steps from her past experiences of identical problems or situations by comparing similarity, relevance, etc. with the experiences \cite{35}. It is shown that the process of solving problems of algebra, i.e., noninsight problems, shows a more incremental pattern than those of insight problems \cite{36}. Thus, we define the following definitions:\\

\begin{definition}
A sequence of thoughts $\{x_i\}_{i=1}^{\infty}$ in a cognitive-consequence space $\mathcal{C}=(C,\sigma, I, Cn)$ is said to converge to a thought $x$ if for  each $\epsilon \in (0,1)$, there exist a positive integer   $m$ such that $Cog(x,x_n)<\epsilon$ for all $n \geq m$. In this case, $x$ is called a cognitive limit of $\{x_i\}_{i=1}^{\infty}$.
\end{definition}

\begin{definition}
Let $\{x_i\}_{i=1}^{\infty}$ be a sequence of thoughts  $\mathcal{C}=(C,\sigma, I, Cn)$ and    $x$ be any thought in a cognitive-consequence space. Then, $x$ is said to be cognitive limit point of $\{x_i\}_{i=1}^{\infty}$ if $B(x, \epsilon)$ contains infinitely many thoughts of $\{x_i\}_{i=1}^{\infty}$ for $\epsilon \in (0,1)$. 
\end{definition} 

If there does not exist such positive integer $m$, then we say that the sequence of thoughts does not converge to its cognitive limit. In that case, there may be a conclusion  that the person fails understand or think and thus fails to attain the cognitive limit. In mathematics, the concept of a Cauchy sequence is present in the study of sequences of real numbers. But the concept of Cauchy sequence is not relevant here since we assign weights to the thoughts and we cannot put thoughts in the real line.\\

\begin{theorem}
If a sequence of thoughts in a cognitive-consequence space $\mathcal{C}=(C,\sigma, I, Cn)$ converges to two cognitive limits  $x'$ and $x''$, then $x'$ and $x''$ cognitively coincide.
\end{theorem}
\begin{proof}
Consider a convergent sequence of thoughts $\{x_i\}_{i=1}^{\infty}$ having two cognitive limits $x'$ and $x''$.
Then for  $\epsilon \in (0,1)$, there exists  a positive integer $N$ such that $Cog(x',x_n)<{\epsilon/2}$ for all $n\geq N'$, and there exists a positive integer $N''$ such that $Cog(x'',x_n)<{\epsilon/2}$ for all $n\geq N''$. We consider $N=max\{N', N''\}$. Then, for $n\geq N$, we apply triangle inequality. 
\begin{align*}
    Cog(x',x'') & \leq Cog(x',x_n)+Cog(x_n,x'')\\
     & = Cog(x',x_n)+Cog(x'',x_n) \\
     & < {\epsilon/2}+{\epsilon/2} \\
     & = \epsilon
\end{align*}

\noindent Since $\epsilon$ is arbitrary, $Cog(x', x'')=0 \Rightarrow x' \approx x''$. Thus, $x'$ and $x''$ cognitively coincide.
\end{proof}

\begin{figure}[h]
    \centering
    \includegraphics[width=14cm, height=11cm]{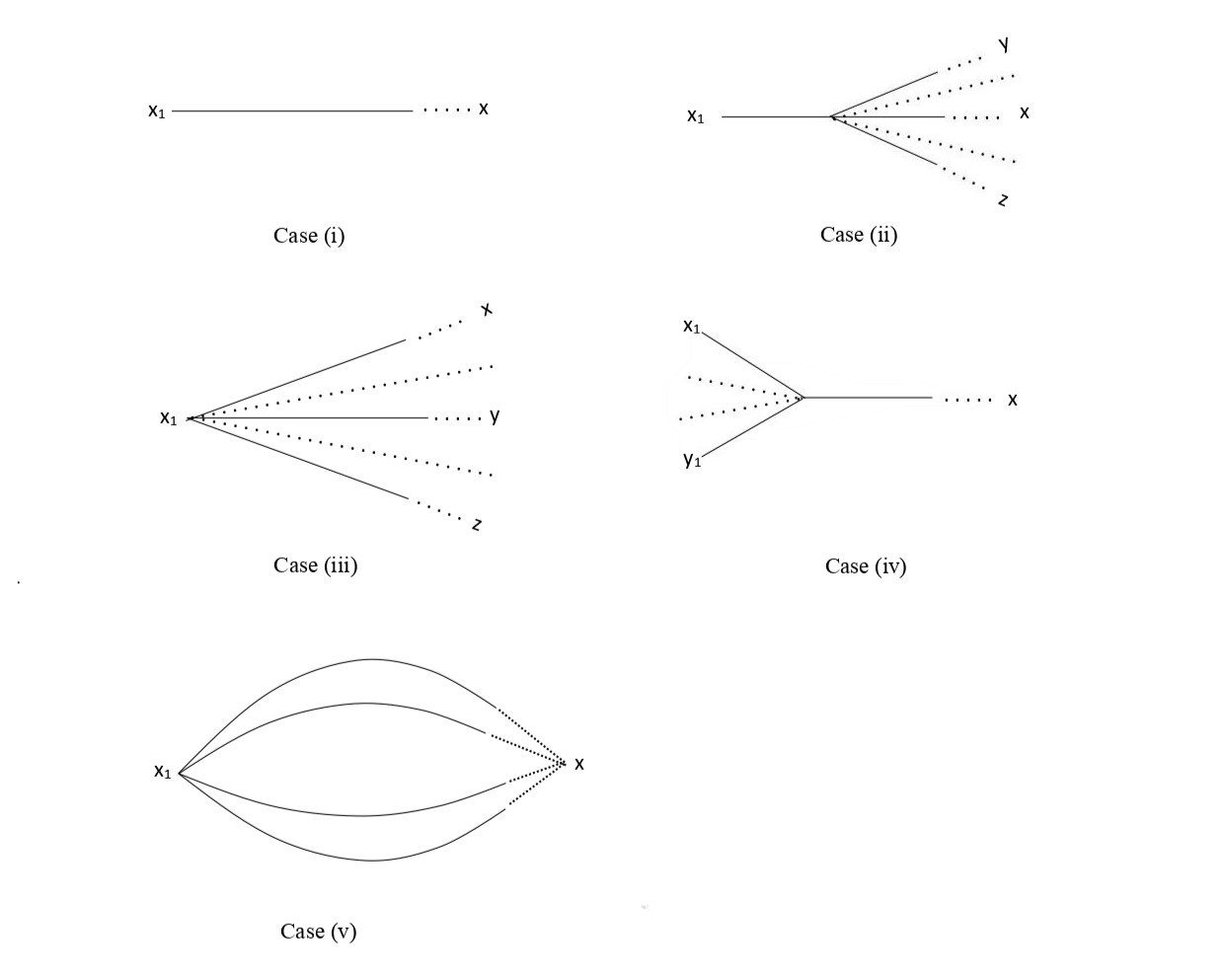}
    \caption{(i) A sequence of thoughts converging to a unique cognitive limit, (ii) a sequence of thoughts eventually diverges to different sequences of thoughts each converging to different or same cognitive limits, (iii) more than one sequences of thoughts diverge from an initial thought and converge to different or same cognitive limit, (iv) more than one different sequences of thoughts merge together to converge to a cognitive limit, (v) more than one sequences of thoughts starting from an initial thought and eventually converge to a  cognitive limit.}
    \label{fig1}
\end{figure}

\noindent \textcolor{red}{From the above theorem, it is clear that the cognitive limit of a sequence of thoughts may not be unique, like the limit of a sequence of real numbers. Due to being semantically similar \cite{18}, cognitively similar \cite{19}, or identical \cite{19}, all the cognitive limits of a sequence of thoughts matter. For example, if a sequence of thoughts converges to two thoughts, say `home' and `house', then it is clear to find that home $\approx$ house, but the meanings of home and house are not identical. Thus, it is important to consider all the cognitive limits of a sequence of thoughts if semantically similar or cognitively similar cognitive limits exist.} According to Rips \cite{37}, in order to arrive at a necessarily true conclusion of an argument, one must be able to construct a mental proof of the conclusion. One cannot construct the mental proof without  prior knowledge of the inference rules needed to complete the proof.  However, a thought may be the cognitive limit point of more than one sequence of thoughts in a cognitive-consequence space. The process of arriving at a conclusion can be applied to the process of arriving at the cognitive limit of a sequence of thoughts. We come across different types of scenarios while arriving at the cognitive limit, which are shown in figure \ref{fig1}. We give the following examples for a clear understanding of each of the cases: \\

\noindent Case (i) is the simplest of all. This type of sequence of thoughts starts with a single thought $x_1$ and finally arrives at the cognitive limit $x$. For example, identifying a digit, doing a simple task such as lifting a pen, book, etc.\\

In case (ii), the sequence initially starts with a single thought, $x_1$, and after some thoughts, it eventually diverges to different sequences of thoughts, each converging to different or the same cognitive limits, $x$, $y$, $z$, etc. For example, to find the solution to a given system of linear equations, we start with an initial thought, and then, after some steps, we can proceed to solve the problem using various methods like graphical methods, substitution, elimination, etc. Another example of this type of sequence is when a student wants to study a specific course in order to get his undergraduate degree. While doing so, the beginning thought works in a single sequence, starting with the initial thought, say, $x_1= ${\it ``I want to study the course A"}. In the next few thoughts, the student will think about choosing an institution to which she wants admission. The sequence will be diverged with the choice of the institution, which leads to different sequences until she gets the undergraduate degree. \\

Case (iii) is somewhat similar to the above case. The only difference is that in this scenario, the sequence of thoughts diverges immediately after the initial thought. $x_1$. For example, to find the initial basic feasible solution of a transportation problem, we start with the given problem and then resort to different methods such as the North West Corner Method, the Least Cost Cell Method, Vogel’s Approximation Method, etc., which produce different initial basic feasible solutions. A hungry tourist in the UK standing in front of restaurants, deciding what to eat, is a non-mathematical example of this kind of scenario. It starts with an initial thought, $x_1=${\it``I am hungry"}, and then depending on her choice of cuisine, say, Indian, Chinese, Italian, etc., she will pick her destination to eat, which gives different sequences of thoughts that start with the thought $x_1$. \\

In case (iv), more than one sequence of thoughts converges into one single sequence, which eventually leads to a single cognitive limit. For example, if we have to solve a mathematical problem that requires knowledge of more than one area, then this kind of convergence takes place. Again, let us give a real-life example, such as detecting the underlying meaning of a painting in an art exhibition. There will be different sequences of thoughts for each element of the painting that eventually merge together in order to detect the overall meaning of the painting. \\ 

The case (v) represents the scenario where the sequence starts with an initial thought $x_1$ and then follows different sequences where each sequence ultimately converges to the same cognitive limit $x$. A mathematical example of this case is to optimise the objective function of a given linear programming problem. We start with the problem, then apply a method of our choice, like the graphical method, the simplex method, etc., and then arrive at the optimal solution. Another example is to find the solution to a simple equation, say, $x+2=3$. Again, suppose a person wants to travel to a specific city from his hometown. Then, the person can take different routes based on feasibility. Now, we have the following theorems:\\

\begin{theorem} \label{t8}
Let $\{x_i\}_{i=1}^{\infty}$ be a sequence of thoughts in a cognitive-consequence space $\mathcal{C}=(C,\sigma, I, Cn)$. If $Cn(\{x_i\}_{i=1}^{\infty})=\{x_i\}_{i=1}^{\infty}$, then the cognitive limit of $\{x_i\}_{i=1}^{\infty}$ belongs to the sequence of thoughts $\{x_i\}_{i=1}^{\infty}$.
\end{theorem}

\begin{proof}
Let $\{x_i\}_{i=1}^{\infty}$ be a sequence of thoughts in a cognitive-consequence space $\mathcal{C}=(C,\sigma, I, Cn)$. \\
Suppose, $Cn(\{x_i\}_{i=1}^{\infty})=\{x_i\}_{i=1}^{\infty}$. To show that the cognitive limit  of $\{x_i\}_{i=1}^{\infty}$ belongs to the sequence of thoughts $\{x_i\}_{i=1}^{\infty}$. That is, if $x$ is the cognitive limit  of $\{x_i\}_{i=1}^{\infty}$, then $x\in \{x_i\}_{i=1}^{\infty}$.\\
Suppose, $x \notin \{x_i\}_{i=1}^{\infty}\Rightarrow x \notin Cn(\{x_i\}_{i=1}^{\infty})$. That is, no previous thought or sequence of thoughts gives the  $x$.  But it is a contradiction as $x$ is the cognitive limit  of $\{x_i\}_{i=1}^{\infty}$. Hence, $x\in \{x_i\}_{i=1}^{\infty}$.
\end{proof}


\begin{theorem} \label{3.10}
Let $\{x_i\}_{i=1}^{\infty}$ be a sequence of thoughts in a cognitive-consequence space $\mathcal{C}=(C,\sigma, I, Cn)$, and let $D$ be the set of cognitive limit(s) of $\{x_i\}_{i=1}^{\infty}$. Then, $Cl^\Box(\{x_i\}_{i=1}^{\infty}) = (\{x_i\}_{i=1}^{\infty}) \bigcup D$.
\end{theorem}
\begin{proof}
Let $\{x_i\}_{i=1}^{\infty}$ be a sequence of thoughts in a cognitive-consequence space $\mathcal{C}=(C,\sigma, I, Cn)$ and let $D$ be the set of cognitive limit(s)  of the sequence. Now, the cognitive closure of the sequence of thoughts $Cl^\Box(\{x_i\}_{i=1}^{\infty})$ is the smallest deductive system which contains $\{x_i\}_{i=1}^{\infty}$ and, from corollary \ref{c3}, $Cl^\Box(\{x_i\}_{i=1}^{\infty})=Cn(\{x_i\}_{i=1}^{\infty})$.\\
From the properties of the consequence operator, we get, $\{x_i\}_{i=1}^{\infty} \subseteq Cn(\{x_i\}_{i=1}^{\infty})$. Since cognitive limit(s) of a sequence of thoughts can be obtained from the consequence of the train of thoughts, so, we can write $D \subseteq Cn(\{x_i\}_{i=1}^{\infty})$. Thus, we get 
\begin{equation} \tag{1}
   \{x_i\}_{i=1}^{\infty} \bigcup D \subseteq Cn(\{x_i\}_{i=1}^{\infty})
\end{equation}
Again, $\{x_i\}_{i=1}^{\infty} \subseteq \{x_i\}_{i=1}^{\infty} \bigcup D \Rightarrow Cn(\{x_i\}_{i=1}^{\infty}) \subseteq Cn(\{x_i\}_{i=1}^{\infty} \bigcup D)$.\\
Since a cognitive limit is inferred from a sequence of thoughts  through some logical reasoning \cite{11}, we get $\{x_i\}_{i=1}^{\infty} \bigcup D$ is a deductive system. Thus, we get
\begin{equation} \tag{2}
    Cn(\{x_i\}_{i=1}^{\infty}) \subseteq Cn(\{x_i\}_{i=1}^{\infty} \bigcup D)= \{x_i\}_{i=1}^{\infty} \bigcup D
\end{equation}
From (1) and (2), we get $Cn(\{x_i\}_{i=1}^{\infty}) = \{x_i\}_{i=1}^{\infty} \bigcup D$. Hence, $Cl^\Box(\{x_i\}_{i=1}^{\infty})= \{x_i\}_{i=1}^{\infty} \bigcup D$.
\end{proof}

\begin{theorem} \label{3.11}
Let $A$ be a deductive system in $\mathcal{C}=(C,\sigma, I, Cn)$ and $\{x_i\}_{i=1}^{\infty}$ be a sequence of thoughts in $A$. Then, the cognitive limit of the sequence of thoughts lies in $A$.
\end{theorem}

\begin{proof}
Let $A$ be a deductive system. So, $Cn(A)=A$. Now, we have the following two cases: \\
{\bf Case (i)}: If the sequence of thoughts $\{x_i\}_{i=1}^{\infty}$ is itself a deductive system, then due to theorem \ref{t8}, the cognitive limit belongs to $\{x_i\}_{i=1}^{\infty}$ as well as the deductive system $A$.\\
{\bf Case (ii)}: Let the sequence of thoughts is non-deductive. Then, $Cn(\{x_i\}_{i=1}^{\infty}) \neq \{x_i\}_{i=1}^{\infty}$. Again, $\{x_i\}_{i=1}^{\infty} \subseteq A \Rightarrow Cn(\{x_i\}_{i=1}^{\infty}) \subseteq Cn(A) =A$. From theorem \ref{3.10}, $Cn(\{x_i\}_{i=1}^{\infty}) = Cl^\Box(\{x_i\}_{i=1}^{\infty})= (\{x_i\}_{i=1}^{\infty}) \bigcup D$, where $D$ is the set of cognitive limit(s). From this we get, $(\{x_i\}_{i=1}^{\infty}) \bigcup D \subseteq A \Rightarrow D \subseteq A$. Thus, the cognitive limit belongs to the deductive system $A$. 
\end{proof}

\begin{theorem}
If $\{x_i\}_{i=1}^{\infty}$ and $\{y_i\}_{i=1}^{\infty}$ be two different sequences of thoughts in $\mathcal{C}=(C,\sigma, I, Cn)$ and they cognitively coincide after a certain position. If they converge to $x$ and $y$ respectively, then  $x \approx y$.
\end{theorem}
\begin{proof}
Let us consider two sequences of thoughts $\{x_i\}_{i=1}^{\infty}$ and $\{y_i\}_{i=1}^{\infty}$ converging to two cognitive limits $x$ and $y$ respectively, where $x_i \approx y_i$ for some $i=k, k+1, ..., \infty$, and $k \in \mathbb{N}$. For  $\epsilon \in (0,1)$,  there exists a positive integer  $N'$ such that $Cog(x, x_i)<\epsilon/2$ for all $i\geq N'$ and there exists another positive integer $N''$ such that $Cog(y, y_i)< \epsilon/2$ for all $i\geq N''$.\\
If $k \leq N', N''$ or $N' \leq k \leq N''$ or $N'' \leq k \leq N'$, then we take $N=max\{N', N''\}$ and if $k \geq N', N''$, then we take $N=k$ and apply triangle inequality for $i\geq N$.
\begin{align*}
    Cog(x, y) & \leq Cog(x, x_i)+Cog(x_i, y)\\
     & = Cog(x, x_i)+Cog(y_i, y)\\
     & = Cog(x, x_i)+Cog(y, y_i)\\
     & < \epsilon/2 + \epsilon/2 \\
     & = \epsilon
\end{align*}
Since $\epsilon$ is arbitrary, $Cog(x, y)=0 \Rightarrow x \approx y$.
\end{proof}

Here, we do not say that the cognitive limits are unique, but rather that the two cognitive limits are semantically similar, cognitively similar, or identical. For instance, without considering the official or legal issue related to a driving licence, one person who can drive a car in India can
also drive a car in the USA. Here, the task of driving the car is cognitively similar but not equal since she must follow the left side in India while she must follow the right side while driving in the USA.\\

\begin{corollary}
If two sequences of thoughts converging to two cognitive limits, respectively, have common elements up to a fixed position at the beginning and have common elements from a fixed position at the end, then the cognitive limit is semantically similar, cognitively similar, or identical.
\end{corollary}
\begin{proof}
Let the two sequences of thoughts have common elements at the beginning up to a certain position, and after diverging to two different paths, then again coincide from some fixed position. The result can be proved by neglecting the common elements at the beginning and then proceeding as per the above theorem.
\end{proof}

\subsection{Function on the cognitive-consequence space and $R^3(t)$}

Sims et al. \cite{8} discussed that there exist continuous maps from the cognitive space to $R^3(t)$. The main motivation for mapping the cognitive space with $R^3(t)$ is that the human mind executes its thoughts in $R^3(t)$. At the same time, the mind also pays attention to what is happening in $R^3(t)$ resulting in another map from  $R^3(t)$ to the cognitive space, which results in another mapping from $R^3(t)$ to the cognitive space $C$ through means of perception. Mathematically, we formulate this with the help of a function as $f$:$\ C \rightarrow R^3(t)$, where $f$ can be defined by physical actions or verbal or written directions that carry out a mental plan in the real world. For example, suppose one has to cook a dish. Before starting the cooking process in $R^3(t)$, she will think of the ingredients and the procedure he will apply in order to accomplish the task, and accordingly, she will execute the steps in $R^3(t)$. At the same time, while executing the steps of the procedure, the mind will also produce some cognitive thoughts depending on the outcome of what the person is doing in real life. \\ 

Lewin's  behavioural  equation is very significant in studying the connection between the environment and the cognitive space. Lewin \cite{4} stated {\it ``Every psychological event depends upon the state of the person and at the same time on the environment, although their relative importance is different in different cases. Thus, we can state our formula B = f(S) for every psychological event as B = f(P,E). The experimental work of recent years shows more and more this twofold relationship in all fields of psychology. Every scientific psychology must take into account whole situations, i.e., the state of both person and environment. This implies that it is necessary to find methods of representing person and environment in common terms as parts of one situation. We have no expression in psychology that includes both."}\\ 

Qi \cite{38} described the link between a person $P$ and her surrounding environment $E$ in Lewin's behavioral equation as {\it ``In this equation P and E are not independent variables. The structure of the environment and the constellation of forces in it vary with the desires and needs, or in general with the state of the person."}  These changes due to continuity  are spatial or temporal in the maps defined between the cognitive space and $R^3(t)$. Thus,  procure the notions of Sims et al. \cite{8} for the following definition:\\

\begin{definition}
Let $\mathcal{C}=(C,\sigma, I, Cn)$ be a cognitive-consequence space, $E$ be an environment in $R^3(t)$ and $(E,T)$ be the practical topological space. Then, the function $f$:$\ C \rightarrow E$ is cognitive-continuous if for an open set $B \in E$, where $B= \bigcup B_i, B_i$ is complete or a connected practical whole in E, we have $f^{-1}(B)=A$, where A is a CWO set in {\it CCT} $\tau$.
\end{definition}

\begin{theorem}
Let $\mathcal{C}=(C,\sigma, I, Cn)$ be a cognitive-consequence space, $E$ be an environment in  $R^3(t)$ and $(E, T)$ be the practical topological space. If $f$:$\ C \rightarrow E$ be a cognitive-continuous function, then $\{f^{-1}(U)$:$\ U \in T\}$ forms a clopen topology on $C$.
\end{theorem}

\begin{proof}
From proposition \ref{p2.1}, the practical topological space $(E, T)$ is a clopen topological space. Now, the topology generated by the family $\{f^{-1}(U)$:$\ U \in T\}$ is the weak topology $T'$ where $f^{-1}(U)$ is a CWO set in {\it CCT} \cite{39}. Again, the inverse image of each clopen set in $T$ is a clopen set in $T'$. Thus, the topology generated by the inverse image of $f$ is clopen.
\end{proof}

\section{Cognitive filter and cognitive ideal}
\begin{figure}[h]
    \centering
    \includegraphics[width=10cm, height=5.5cm]{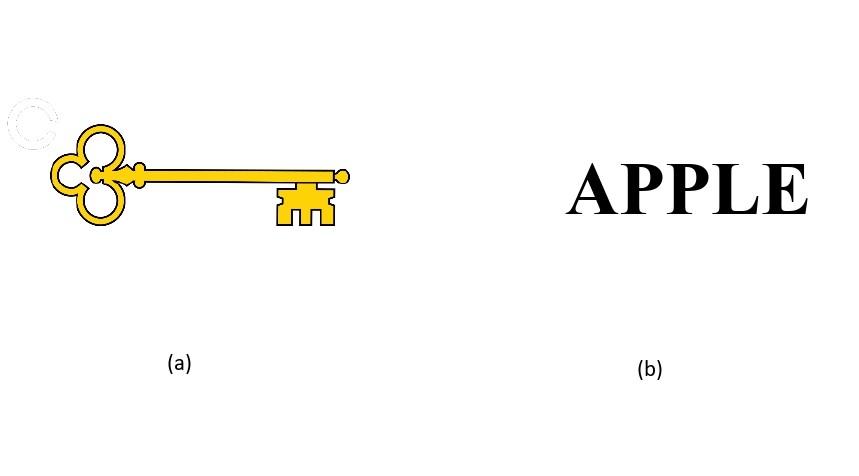}
    \caption{(a) a mental representation of a key, (b) a word say `APPLE'.}
    \label{fig2}
\end{figure}

In order to deal with a real-life problem, one must organise the mental representations related to the problem into some subsets of $\mathcal{C}=(C,\sigma, I, Cn)$ for the easy execution of the solution. These subsets may constitute filters, ideals, etc., or sometimes just a collection of mental representations from the cognitive-consequence space. These subsets can be created based on the problem solver's needs. In this section, we discuss some filters and ideals defined on a cognitive-consequence space $\mathcal{C}=(C,\sigma, I, Cn)$. There may be direct or indirect connections between two mental representations. The direct connections are easy to comprehend. For indirect connections, let us take an example of the word `keyword'. When we think of the words `key' and `word' separately, two totally different mental representations emerge in the mind. But when we combine the two words together, it gives a completely new mental representation of something having an entirely different meaning. Figure \ref{fig2} gives a clear understanding of the above example. Thus, we define the following definition: \\

\begin{definition}
Let $\mathcal{C}=(C,\sigma, I, Cn)$ be a cognitive-consequence space. A subset $A \subseteq C$ is called psychologically connected to the mental representation $f_* \in C$ such that $A$ is a collection of mental representations $\{f_{*_n}\}$ which are connected to the thought $f_*$ through a sequence of thoughts in a direct or indirect manner so that the mental representation $f_*$ makes sense in the real world or in the mind.
\end{definition}

\begin{theorem} \label{4.1}
Let $\mathcal{C}=(C,\sigma, I, Cn)$ be a cognitive-consequence space. If the set $\hat{f_*}$ consists of all subsets $A \subseteq C$ such that $A$ is psychologically connected to a mental representation $f_* \in C$, then the collection $\hat{f_*}$ forms an ideal in the cognitive-consequence space $\mathcal{C}$.  
\end{theorem}

\begin{proof}
The set $\hat{f_*}$ is non-empty since the thought ``$f:$ no thoughts about $f_*$" is also a thought connected to $f_*$. Let $A \in \hat{f_*}$ and $B \subseteq C$. If $B \subset A$, then the mental representations in $B$ are also psychologically connected to $f_*$. Thus, $B \in \hat{f_*}$.\\ \\
Again, let $A, B \in \hat{f_*}$. Since $A \bigcup B$ consists of thoughts that are psychologically connected to $f_*$, thus $A \bigcup B \in \hat{f_*}$.
\end{proof}

\begin{definition}
Let $\mathcal{C}=(C,\sigma, I, Cn)$ be a cognitive-consequence space. A consequence-ideal $\mathcal{I}$ on $\mathcal{C}$ is a collection of subsets of $C$ such that:
\begin{enumerate} [label=(\roman*)]
    \item $Cn(\phi) \in \mathcal{I}$,
    \item  for $A, B \subseteq C$ with $Cn(A) \in \mathcal{I}$ and $Cn(B) \subseteq C$; if $Cn(B) \subset Cn(A)$, then $Cn(B) \in \mathcal{I}$,
    \item if $Cn(A) \in \mathcal{I}$ and $Cn(B) \in \mathcal{I}$, then $Cn(A) \bigcup Cn(B) \in \mathcal{I}$.
\end{enumerate}
\end{definition}
Though the existing definition of ideal \cite{40} and the above-defined consequence ideal are found to be similar at first glance, there is a significant difference between these two definitions. The empty set $\phi$ should be there in the ideal, whereas in the consequence ideal $Cn(\phi)$ belongs to it, which is not equal to $\phi$ but represents the smallest system consisting of all logically true statements since any logically true statement can be derived from an empty set of assumptions \cite{11}. Moreover, $Cn(\phi)$ is the smallest system since $\phi \subseteq A$ for all $A$ in $C$, and from the properties of the consequence operator we get, $Cn(\phi)\subseteq Cn(A)$ for all systems $Cn(A)$ in $C$.\\

\begin{theorem} \label{4.2}
Let $M$ be a collection of meaningful mental representations in a cognitive-consequence space $\mathcal{C} = (C,\sigma, I, Cn)$. Let the set $\Hat{\Hat{f_*}}$ consist of all subsets $A \subseteq M$ such that $A$ is the collection of logically true meaningful mental representations that are connected to $f_* \in M$. Then, the collection $\Hat{\Hat{f_*}}$ forms a consequence-ideal in the cognitive-consequence space $\mathcal{C}$. 
\end{theorem}

\begin{proof}
Logical inferences cannot produce false conclusions from true assumptions \cite{11}. Thus, if $A \in \Hat{\Hat{f_*}}$, then $Cn(A) \in \Hat{\Hat{f_*}}$. Now, $Cn(\phi) \in \Hat{\Hat{f_*}}$ since $Cn(\phi)$ is the system of all logically true statements, i.e., $Cn(\phi)=$ {\bf LOG}. \\
Again, let $Cn(A) \in \Hat{\Hat{f_*}}$ and $Cn(B) \subseteq M$. If $Cn(B) \subset Cn(A)$, then the mental representations in $Cn(B)$ are also true statements connected to $f_*$. Thus, $Cn(B) \in \Hat{\Hat{f_*}}$.\\
Again, let $Cn(A)$, $Cn(B) \in \Hat{\Hat{f_*}}$, then the union $Cn(A) \bigcup Cn(B) \in \Hat{\Hat{f_*}}$ since $Cn(A) \bigcup Cn(B)$ contains logically true statements.
\end{proof}

\begin{definition}
Let $\mathcal{C}=(C,\sigma, I, Cn)$ be a cognitive-consequence space. A consequence filter $\mathcal{F}$ on $\mathcal{C}$ is a collection of subsets of $C$ satisfying the following conditions: 
\begin{enumerate} [label= (\roman*)]
    \item $Cn(C) \in \mathcal{F}$,
    \item for $Cn(A) \in \mathcal{F}$ and $Cn(B) \subseteq C$, if $Cn(A) \subset Cn(B)$, then $Cn(B) \in \mathcal{F}$,
    \item if $Cn(A) \in \mathcal{F}$ and $Cn(B) \in \mathcal{F}$, then $Cn(A) \bigcap Cn(B) \in \mathcal{F}$.
\end{enumerate}
\end{definition}

\begin{theorem} \label{t4.3}
Let $C_d \subset C$ be the deductive part of a cognitive-consequence space $\mathcal{C}=(C,\sigma, I, Cn)$, where $C$ consists of both deductive and non-deductive parts, and let $f \in C_d$ be a mental representation. We define $f_d = \{ A \subseteq C_d: f \in A$, and $A$ contains at least one subset $B$ of $C$ such that $Cn(B)=B\}$. Then, $f_d$ forms a filter on $C_d$.
\end{theorem}

\begin{proof}
$C_d \in f_d $ since $f\in C_d$ and $C_d$ itself is the deductive system. Again, let $A \in f_d$. Then, $f \in A$ and there exist at least one subset $D \subseteq A$ such that $Cn(D)=D$. Let $B \subseteq C_d$ and $A \subset B$. Then, $f \in B$ and there exists $D \subseteq A \subseteq B$ such that $Cn(D)=D$. Hence, $B \in f_d$. \\ \\
Again, let $A, B \in f_d$. Then, $f \in A, B$ and there exist $D \subseteq A$ and $E \subseteq B$ such that $Cn(D)=D$ and $Cn(E)=E$. Now, $D \bigcap E$ is a deductive system since $Cn(D \bigcap E)=Cn(D) \bigcap Cn(E) = D \bigcap E$. Then, $f \in A \bigcap B$ and there exist $D \bigcap E \subseteq A\bigcap B$ such that $Cn(D \bigcap E)=D\bigcap E$. So, $A \bigcap B \in f_d$. Thus, $f_d$ forms a filter on $C_d$.
\end{proof}

\begin{corollary}
The filter defined in theorem \ref{t4.3} is a consequence filter in $C_d$.
\end{corollary}
\begin{proof}
$Cn(C_d) \in f_d$, since $f \in C_d = Cn(C_d)$ and $C_d$ is the deductive system. \\ \\
Again, let $Cn(A) \in f_d$, $Cn(B) \subseteq C_d$ and $Cn(A) \subset Cn(B)$. Then, $f \in Cn(B)$ and $Cn(B)$ itself is a deductive system. So, $Cn(B) \in f_d$.\\ \\
Now, let $Cn(A), Cn(B) \in f_d$. Then, $f \in Cn(A) \bigcap Cn(B)$ and $Cn(A) \bigcap Cn(B)$ is a deductive system. Thus, $Cn(A) \bigcap Cn(B) \in f_d$. Thus, $f_d$ forms a consequence filter on $C_d$.
\end{proof}

\noindent In a similar manner, we may also have such a filter associated with more than one mental representation in $C_d$. For example, $(fg)_d = \{ A \subseteq C_d: f,g \in A$ and $A$ contains at least one subset $B$ such that $Cn(B)=B\}$, $(fgh)_d = \{ A \subseteq C_d: f,g,h \in A$ and $A$ contains at least one subset $B$ such that $Cn(B)=B\}$, etc. The concepts of filter and consequence filter are found to be similar at first glance, but there is a significant difference between these two. It is easy to check that the consequence filter and filter are different from each other.\\

\noindent There are many problems that can be solved using deductive inferences. Given a mental representation $f$, the filter $f_d$ constitutes those subsets $A$ of $C_d$ that contain a deductive system related to $f$. Here, $f$ itself can be inside the deductive system or it may belong to a superset of the deductive system. The main purpose of the consequence filter is to organise the available information in a manner such that it will be useful to solve the problem. This type of organisation of mental representations in terms of a deductive system can be useful in logical decision-making since logical inferences do not produce false conclusions from true assumptions during the development of a theory because a theory is a deductive system of an ordered set of sentences \cite{11}. Since human beings do not always obey deductive rules \cite{11, 17}, defining a consequence filter only in terms of deductive systems may be problematic. Thus, we extend the idea of organising subsets containing mental representations to the whole cognitive-consequence space, irrespective of whether a set of mental representations forms a deductive system or not. Hence, we have the following theorem:\\

\begin{theorem} \label{4.5}
Let $\mathcal{C}=(C,\sigma, I, Cn)$ be a cognitive-consequence space and $f \in C$. Consider the collection $\hat{f} = \{ A \subseteq C : f \in Cn(A) \}$. Then, $\hat{f}$ forms a filter in $C$.
\end{theorem} 
\begin{proof}
$C \in \hat{f}$ since $f \in C \subseteq Cn(C)$.
Again, let $A \in \hat{f}$. Then, $f \in Cn(A)$. We consider $B \subseteq C$ and $A \subset B$. Then, $A \subset B \Rightarrow Cn(A) \subset Cn(B)$. Hence, $f \in Cn(B) \Rightarrow B \in \hat{f}$.
\\ \\
Now, let $A,\ B \in \hat{f}$. Thus, $f \in Cn(A)$ and $f \in Cn(B)$. Hence, $f \in Cn(A) \cap Cn(B)$. So, $f \in Cn(A \bigcap B)$. Hence, $A \bigcap B \in \hat{f}$. Thus, $\hat{f}$ is a filter on $C$.
\end{proof}
The main advantage of defining a filter like this is that we can organise the mental representations associated with a thought $f$ that can be considered both a deductive system and a non-deductive system.\\

\section{G\"{o}del's incompleteness black hole}
In general topology, compactness deals with the concepts of open covers and open subcovers \cite{41}. But in a cognitive sense, we can think of the compactness from the point of view of problem solving. To solve a problem, one must proceed through some steps that eventually lead to the solution of the problem. In most cases, while solving a problem, the solver depends on past experiences that bear a resemblance to the given problem \cite{22}. These past experiences, to some extent, constitute the steps for the solution, which cover the solution space. But, sometimes, after some steps, the problem solver must halt since her next step does not match any prior experiences or she becomes clueless about what to do next. We have the following real-life situation to understand: \\ \\
Fermat's Last Theorem, one of the most notable and hardest theorems in mathematical history, which was unsolved for more than three hundred years, became provable after Gerhard Frey claimed that the proof of Fermat's Last Theorem was the direct consequence of the proof of the Taniyama-Shimura conjecture \cite{42}. To put it simply, the idea was that if one considered Fermat's Last Theorem to be false, then the Taniyama-Shimura conjecture would also be false. Equivalently, if the Taniyama-Shimura conjecture could be proven to be true, then Fermat's Last Theorem would be true. Following this remarkable claim, the unsolved proof became plausible. In this famous example, we are trying to convey that sometimes an unsolvable problem can be solved if one can find an equivalent way to solve the problem.\\

\begin{figure}[h]
    \centering
    \includegraphics[width=10cm, height=6cm]{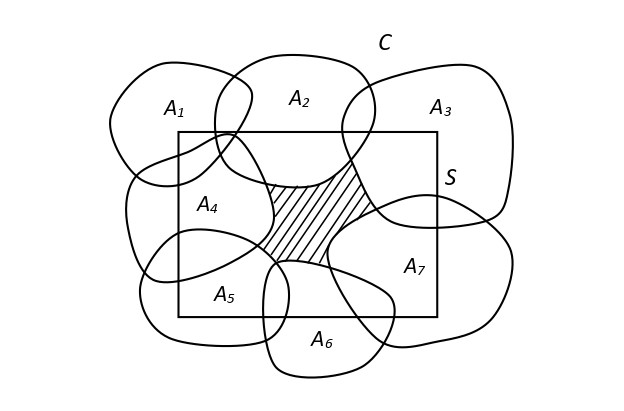}
    \caption{G\"{o}del's incompleteness black hole}
    \label{fig3}
\end{figure}

\noindent \\ Now, we are going to show the process of problem solving using one simple example of having one equivalent way of finding the solution. Mathematically, we can represent the process of problem solving as a sequence of mental representations $\{x_i\}_{i=1}^{\infty}$, where each $x_i \in C$ is a step to solve the problem. In figure \ref{fig3}, we consider a solution space $S$ of a problem $P$ in a cognitive-consequence space $\mathcal{C}=(C,\sigma, I, Cn)$. Let $\{A_i|\; A_i \subseteq C,\ i \in \Delta,\; \Delta$ is an index set\} be a cover of $S$ which contains a sequence $\{x_i\}_{i=1}^{\infty}$ of mental representations of the steps of the solution of a problem $P$. Let the finite sub-collection $\{A_1, A_2, \dots, A_n\}$ contains the known mental representations from one's past experiences that are used to solve the problem. Then, they contain a finite number of thoughts from the sequence, i.e., $\{x_i\}_{i=1}^{m} \in \bigcup_{i=1}^{n} \{A_i\}$. Now, $\{x_i\}_{i={m+1}}^{\infty} \in \bigcup_{i=n+1}^{\infty} \{A_i\}$ is completely unknown to the solver. In order to solve the problem, the solver can find another known cover $\{B_j|\; B_j \subseteq C,\ j \in \Delta',\; \Delta'$ is an index set\} of $S$ that contains a sequence of thoughts $\{y_j\}_{j=1}^{\infty}$ giving the solution to the problem in an alternate way. Then, we say this cover cognitively coincides with the cover $\{A_i|\; A_i \subseteq C,\ i \in \Delta,\; \Delta$ is an index set\} and can be written as $\bigcup_{j=1}^{\infty} \{B_i\} \approx \bigcup_{i=1}^{\infty} \{A_i\}$. The thoughts in $\{y_j\}_{j=1}^{\infty}$ and $\{x_i\}_{i={1}}^{\infty}$ are not necessarily cognitively similar or identical, but what is important is that the cognitive limit in both cases should be cognitively similar. \\

There is a lot of debate going on regarding the validation of G\"{o}del's incompleteness theorem \cite{43} in the case of the human mind. According to the G\"{o}dels incompleteness theorem \cite{43}, there is a statement that is neither provable nor refutable in the given axiomatic system. Lucas \cite{44} stated that {\it ``G\"{o}del’s theorem applies to deductive systems, and human beings are not confined to making only deductive inferences. G\"{o}del’s theorem applies only to consistent systems, and one may have doubts about how far it is permissible to assume that human beings are consistent."} Initially, the shaded region in Figure 3 represents those mental representations of $C$ that are not included in any past experience but exist in the solution space. This shaded region of incomplete information or incomplete steps helps to satisfy G\"{o}del's incompleteness theorem. The existence of the required solution to a given problem is determined by this region. We call this shaded region as G\"{o}dels's incompleteness black hole. Thus, we define the following definition:\\

\begin{definition}
In the solution space $S$ of a problem $P$, a G\"{o}del's incompleteness black hole $A \subset C$ exists if there is a solution sequence of thoughts  $\{x_n\}_{n=1}^{\infty}$ with a virtual cognitive limit $x$, and for an $\epsilon \in (0,1)$, there exists a positive integer $k$ such that $x_n \notin B(x, \epsilon) \subseteq A$ for all $n\geq  k$. \end{definition}

\begin{definition}
The solution space $S$ of a problem $P$  in a cognitive-consequence space $\mathcal{C}=(C, \sigma, I, Cn)$ is said to be cognitively compact if there does not exist any G\"{o}del's imcompleteness black hole in it. 
\end{definition}

\begin{theorem} \label{5.1}
If a G\"{o}del's incompleteness black hole $A$ exists in the solution space $S$ of a problem $P$, then the solution sequence of thoughts is not convergent in $A$. 
\end{theorem}
\begin{proof}
Suppose a  G\"{o}del's incompleteness black hole  $A$ exists in the solution space $S$ of a problem $P$ with a solution sequence of thoughts $\{x_n\}_{n=1}^{\infty}$ and a virtual cognitive limit $x$. Then, for an $\epsilon \in (0,1)$, there exists a positive integer $k$ such that $x_n \notin B(x,\epsilon) \subseteq A$ for all $n \geq k$, i.e.,  $x_n \notin \{y \in C: Cog(x, y) < \epsilon\}$ for all $n \geq k$, i.e., $Cog (x, x_n)\geq \epsilon$ for all $n> k$. Thus, for the above $\epsilon$, one will never reach to the last thought which yields immediately her solution $x$.
\end{proof}

\section{The limit of artificial intelligence}
In artificial intelligence, the term cognitive computing is used to discuss the adaptive nature of artificial intelligence \cite{3}. That is, they reason and learn through their interactions with humans as well as from their experiences in their surroundings, as opposed to being explicitly programmed. In this case, there is no predetermined set of axioms from which all possible inferences may be derived, in contrast to using an algorithm to identify a machine \cite{3}. Although it can be claimed that an artificial intelligence system is adaptive, it may include axiomatic subsystems \cite{2}. G\"{o}del's  incompleteness theorem is applicable to deductive or axiomatic systems, which consist of a few axioms and rules of inference. According to Turing, ``to solve a problem, a computer must follow a predetermined algorithm" \cite{2}. If a machine follows a specified algorithm comprising of a set of instructions or axioms, then it becomes an axiomatic system and so must satisfy G\"{o}del's incompleteness theorem \cite{44}. Unlike human intelligence, this leads to some limitations for artificial intelligence. For instance, in July 2022, a chess-playing robot broke a seven-year-old player's finger during a match at the Moscow Chess Open \cite{45}. Sergey Lazarev, President of the Moscow Chess Federation, stated, ``The child made a move, and after that, we needed to give time for the robot to answer, but the boy hurried and the robot grabbed him.'' This incidence  concludes that algorithmic instructions are not enough in the case of artificial intelligence; also,  human-like consciousness and humanity are required. It is found that human-like consciousness is influenced by the surroundings of a human being. Moreover, religions and religious thoughts play a major role in establishing humanity in a human being, either consciously or subconsciously \cite{46}. Instabilities that arise from deep learning as well as in modern artificial intelligence are a very common thing in today's artificial intelligence world. Colbrook et al. \cite{47} came to the conclusion in 2022 that even though precise neural networks exist, there are some algorithms that cannot compute well-conditioned problems. They demonstrated some fundamental restrictions on the existence of algorithms that can calculate the stable neural networks that are now in existence, thereby limiting artificial intelligence, and gave an answer to  Smale's eighteenth problem \cite{1} in the case of the limit of artificial intelligence. It is also suspected that the fuzzy nature of the available data will make this limit more concrete. Thus, it is understood from the above-stated evidence that there is a limit to artificial intelligence.\\

\section{Discussion}
 Due to Smale \cite{1}, many mathematicians have shown keen interest in continuing research on the limits of human intelligence as well as artificial intelligence. On the other hand, Luwin's \cite{4} explanations to study the cognitive dynamics and related psychological processes of the human mind, by using the ideas of  topology and relativity, have encouraged many experts  to use topology as a useful tool to study cognitive behaviour, human intelligence, and artificial intelligence. But, in reality, it may be a difficult task to know the limit of human intelligence using only experiments. According to Thagard \cite{6}, Von Neumann \cite{17}, Lucas \cite{44}, and many others, a human mind is not a machine. A human mind  does not only follow deductive inferences. Since consequence operators  deal with meaningful sentences \cite{16},  it is important to study the eighteenth problem  of Smale \cite{1} theoretically from the perspective of topology and consequence operators. As discussed in the previous section, Colbrook et al. \cite{47} provided  an answer to Smale's eighteenth problem \cite{1} in the case of the limit of artificial intelligence.  Interestingly, theorems \ref{t3} and \ref{t4}  give answers to the eighteenth problem of Smale \cite{1} in the case of human intelligence. Since at least one mental representation does not exist  in any CWO set, and at least one mental representation always exists in a CWO set, it can be concluded that human intelligence is  limitless. For a layman, it can be stated that at least one metal representation will always remain in the deductive part, and at least one metal representation will always remain in the non-deductive part.  Moreover, Griffiths \cite{48} studied human intelligence from the perspective of three limitations, viz., limited time, limited computation, and limited communication. According to him, a limited human life span is one of the factors contributing to the limitation of human intelligence. Griffiths \cite{48} also pointed out that limited time yields limited amounts of data. Thus, a human being is to perform all her tasks, to which intelligence is to be applied, with limited amounts of data. However, we provide an example  to contradict Griffiths \cite{48}, and it demonstrates the superiority of human intelligence over artificial intelligence. In the era of data-driven artificial intelligence, ancient texts like Surya Siddhanta, the Pur\=a\d nas, etc. bear the evidence of human intelligence,  which almost accurately predicted several astronomical facts without any prior data. For example, in an ancient Indian text by S\={a}ya\d na (c. 1315-1387), the speed of light was calculated as 186,000 miles per second, which is equivalent to 299,337,984 metres per second \cite{49}. This calculated value of the speed of light is very close to the correct value of the speed of  light, i.e.,  $3 \times 10^8$ metres per second. Thus, we suggest looking into the limits of human intelligence in the macro sense in lieu of the  micro sense. For example, in physics, Newton \cite{50} laid the foundations of  Newtonian mechanics. After several scientific contributions by others, Einstein \cite{51} introduced relativity theory  to reveal the various physical phenomena related to the universe. Since then, many researchers have been exploring various new physical phenomena  related to  the universe. It is also an example of the limits of intelligence in   macro sense. This paper provides an answer to the limit of human intelligence \cite{1} in a macro sense. Hence, we also conclude that human intelligence will always remain superior to artificial intelligence in various aspects.\\

On the other hand, a human being is capable of dealing with a situation or a problem in various ways. This inspires us to study the process as a sequence of thoughts occurring sequentially within the mind. While doing this, essential concepts like cognitively similar, semantically similar, etc. come into the picture. Since a human being uses both deductive and non-deductive approaches in different situations, it is also important to organise the mental representations in a way that they can be useful in the preferred approach to problem solving. Theorems \ref{4.1}, \ref{4.2} and \ref{4.5} deal with the non-deductive approach, while theorem \ref{t4.3} is effective to organise the mental representations when one is going to use the deductive approach. Moreover, G\"{o}del's incompleteness black hole is one of the real scenarios that the human mind often faces. Thus, theorem \ref{5.1} connects the problem-solving approach to any problem in case of the non-availability of the proper sequence of thoughts to reach the conclusion. \\

\section{Conclusion}
In this paper, we formulate the mental space of a human being as a cognitive structure $\mathcal{C}=(C, \sigma, I, Cn)$ called a cognitive-consequence space as a generalisation of the mental structure defined by Sims et al. \cite{8}. We consider Tarski's consequence operator and Lewin's topological psychology as useful theoretical foundations to proceed with all the results and justifications of this paper. We construct a cognitive-consequence topological space $(C, \tau)$  and discuss some fundamental properties. We discuss the notion of the cognitive limit of a sequence of thoughts in a cognitive-consequence space and illustrate different scenarios for arriving at that cognitive limit. Further, we study cognitive-continuous functions defined between a cognitive-consequence space and $R^3(t)$ establishing the relationship between the mind and the environment. We discuss the notions of cognitive ideal and cognitive filter and give some examples in the cognitive-consequence space. Lastly, we discuss the existence of G\"{o}del's incompleteness black hole in the solution space of a given problem. However, theorems \ref{t3} and \ref{t4} conclude that human intelligence is limitless, and thus, we provide an answer to the eighteenth problem of Smale \cite{1} in the case of human intelligence.  To the best of our knowledge, no mathematical proof was provided earlier, prior to us working on it. We also provide justifications for the limitations of artificial intelligence. Hence, we hope this paper will be interesting as well as important to experts in many related interdisciplinary areas. \\

{\bf Conflict of interests: } The authors declare that there is no conflict of interest.\\

{\bf Funding: } This research has not received any funding.

\end{document}